\documentclass{amsart}
\usepackage{hyperref}
\usepackage{amsmath, amsthm, amsfonts, mathtools, tikz-cd, amssymb}
\usepackage{enumitem}

\makeatletter
\@namedef{subjclassname@2020}{\textup{2020} Mathematics Subject Classification}
\makeatother

\title{Amenable covers of right-angled Artin groups}
\author{Kevin Li}
\address{School of Mathematical Sciences, University of Southampton, Southampton SO17 1BJ, United Kingdom}
\email{kevin.li@soton.ac.uk}
\date{\today}


\theoremstyle{definition}
\newtheorem{defn}{Definition}[section]

\theoremstyle{plain}
\newtheorem{thm}[defn]{Theorem}
\newtheorem{lem}[defn]{Lemma}
\newtheorem{prop}[defn]{Proposition}
\newtheorem{cor}[defn]{Corollary}
\theoremstyle{remark}
\newtheorem{rem}[defn]{Remark}

\newcommand{\IN}{\ensuremath\mathbb{N}}
\newcommand{\IZ}{\ensuremath\mathbb{Z}}
\newcommand{\IR}{\ensuremath\mathbb{R}}
\newcommand{\IF}{\ensuremath\mathbb{F}}
\newcommand{\F}{\ensuremath\mathcal{F}}
\newcommand{\E}{\ensuremath\mathcal{E}}
\newcommand{\FS}{{\ensuremath\mathcal{F}\spann{\calS}}}
\newcommand{\TR}{\ensuremath\mathcal{TR}}
\newcommand{\FIN}{\ensuremath\mathcal{FIN}}
\newcommand{\AME}{\ensuremath\mathcal{AME}}
\newcommand{\calS}{\ensuremath\mathcal{S}}
\newcommand{\calH}{\ensuremath\mathcal{H}}
\newcommand{\TC}{\ensuremath\mathsf{TC}}
\newcommand{\Subexp}{\ensuremath\mathsf{Subexp}}
\newcommand{\EFG}{\ensuremath E_\mathcal{F}G}

\newcommand{\enum}{\rm{(\roman*)}}
\newcommand{\spann}[1]{{\ensuremath \langle{#1}\rangle}}
\newcommand{\into}{\ensuremath\hookrightarrow}

\DeclareMathOperator{\cd}{cd}
\DeclareMathOperator{\gd}{gd}
\DeclareMathOperator{\cat}{cat}
\DeclareMathOperator{\vcd}{vcd}
\DeclareMathOperator{\id}{id}
\DeclareMathOperator{\LScat}{LS-cat}
\DeclareMathOperator{\im}{im}

\begin{document}

\maketitle

\begin{abstract}
	Let~$A_L$ be the right-angled Artin group associated to a finite flag complex~$L$.
	We show that the amenable category of~$A_L$ equals the virtual cohomological dimension of the right-angled Coxeter group~$W_L$.
	In particular, right-angled Artin groups satisfy a question of Capovilla--L\"oh--Moraschini proposing an inequality between the amenable category and Farber's topological complexity. 
\end{abstract}

\section{Introduction}

A classical approach to study a topological space~$X$ is to cover it by open subsets~$U_0,\ldots,U_n$ that are simpler or small in an appropriate sense and to analyse how these overlap. The minimal possible cardinality~$n$ of such a cover yields a measure of complexity of the space~$X$.
When the subsets~$(U_i)_i$ are required to be contractible in~$X$, we obtain the Lusternik--Schnirelmann category (LS-category for short)~$\LScat(X)$ which is a well-studied homotopy invariant originating from critical point theory~\cite{CLOT03}. 
We will relax the contractibility assumption and instead require the subsets~$(U_i)_i$ to be \emph{amenable in~$X$}, in the sense that the group
\[
	\im\bigl(\pi_1(U_i\into X,x)\bigr)
\]
is amenable for every basepoint~$x\in U_i$. Then the \emph{amenable category} $\cat_\AME(X)$ of~$X$ is the minimal~$n\in \IN_{\ge 0}$ for which there exists an open cover~$X=\bigcup_{i=0}^n U_i$ by~$n+1$ many amenable subsets. 
Clearly, we have~$\cat_\AME(X)\le \LScat(X)$.

Amenable groups (such as finite or abelian groups) and hence amenable subsets can be considered as small for many purposes in geometry, topology, and dynamics.
Therefore the amenable category is a meaningful threshold, especially for aspherical spaces. For instance, there are vanishing results in all degrees larger than the amenable category for the comparison map from bounded cohomology to singular cohomology~\cite{Gromov82,Ivanov85}, for~$\ell^2$-Betti numbers~\cite{Sauer09}, and for homology growth~\cite{Sauer16,Haulmark-Schreve22}.
The amenable category was systematically studied as an invariant of 3-manifolds in~\cite{GGH13,GGH14} and for arbitrary spaces recently in~\cite{CLM20,Loeh-Moraschini21}.

The focus of this note is on the amenable category of aspherical spaces. Since the amenable category is a homotopy invariant, it yields an invariant of discrete groups~$G$ by setting~$\cat_\AME(G)\coloneqq \cat_\AME(BG)$. Here~$BG$ is an Eilenberg--MacLane space. 
By the classical work of~\cite{Eilenberg-Ganea57,Stallings68,Swan69}, the LS-category $\LScat(BG)$ coincides with the cohomological dimension~$\cd(G)$. In particular, we always have~$\cat_\AME(G)\le \cd(G)$.
The amenable category is difficult to compute in general, the usual strategy being to exhibit an explicit open cover by amenable subsets and to prove its minimality using (co)homological obstructions.
The precise value of~$\cat_\AME(G)$ is known, e.g., for the following classes of groups:
\begin{itemize}
	\item $\cat_\AME(G)=0$ if and only if~$G$ is amenable;
	\item $\cat_\AME(G)=1$ if and only if~$G$ is a non-amenable fundamental group of a graph of amenable groups~\cite[Corollary~5.4]{CLM20};
	\item $\cat_\AME(G)=\cd(G)$ if~$G$ is torsion-free non-elementary hyperbolic \cite{Mineyev01} \cite[Example~7.8]{CLM20}. 
\end{itemize}

The main result of the present note is a computation of the amenable category for all right-angled Artin groups. These form an important class of groups in geometric group theory, interpolating between free groups and free abelian groups. Let~$L$ be a finite flag complex (i.e., a simplicial complex in which every clique spans a simplex) with vertex set~$V$. The right-angled Artin group~$A_L$ has as generators vertices~$v\in V$, subject to the relation that~$v_1$ and $v_2$ commute if and only if they are connected by an edge in~$L$. The right-angled Coxeter group~$W_L$ is the quotient of~$A_L$ obtained by adding the relations that each generator~$v\in V$ is of order~2. Since~$W_L$ is virtually torsion-free, its virtual cohomological dimension~$\vcd(W_L)$ is well-defined as the cohomological dimension of a finite index torsion-free subgroup.

\begin{thm}[Corollary~\ref{cor:main}]\label{thm:main intro}
	Let~$A_L$ be the right-angled Artin group associated to a finite flag complex~$L$.
	Then we have
	\[
		\cat_\AME(A_L)=\vcd(W_L) \,.
	\]
\end{thm}

Theorem~\ref{thm:main intro} provides many examples of groups for which the amenable category is not extremal, in the sense that~$1<\cat_\AME(G)<\cd(G)$.
Furthermore, it follows from Theorem~\ref{thm:main intro} and~\cite{Dranishnikov97} that there are right-angled Artin groups~$A_{L_1}$ and~$A_{L_2}$ satisfying~$\cat_\AME(A_{L_1}\times A_{L_2})<\cat_\AME(A_{L_1})+\cat_\AME(A_{L_2})$.

Another invariant of a similar spirit is Farber's topological complexity~$\TC$ which is motivated by the motion planning problem in robotics~\cite{Farber03}.
In~\cite[Question~8.1]{CLM20} it is asked for which topological spaces~$X$ the following inequality holds:
\[
	\cat_\AME(X\times X)\le \TC(X) \,.
\]
Examples of spaces and groups satisfying this inequality can be found in~\cite[Section~8]{CLM20}, and no counter-example seems to be known at the time of writing.
We show that all right-angled Artin groups are positive examples. 

\begin{thm}[Proposition~\ref{prop:question}]
	Let~$A_L$ be the right-angled Artin group associated to a finite flag complex~$L$.
	Then we have~$\cat_\AME(A_L\times A_L)\le \TC(A_L)$.
\end{thm}

We also obtain a complete characterisation of right-angled Artin groups with (non-)vanishing minimal volume entropy (Theorem~\ref{thm:minvolent}),
resolving the cases that were not covered by recent work in~\cite{Haulmark-Schreve22,Bregman-Clay21}.

Our proofs rely on combining upper and lower bounds (Lemma~\ref{lem:basic}) with existing results on generalised LS-category, classifying spaces for families of subgroups, and homology growth from~\cite{CLM20, Haulmark-Schreve22, Loeh-Moraschini21, Okun-Schreve21, Petrosyan-Prytula20, Sauer16}.

\subsection*{Acknowledgements}
The present note is part of the author's PhD project.
He thanks his advisors Ian Leary for extremely helpful discussions and Nansen Petrosyan for his support. 
We are grateful to Mark Grant, Sam Hughes, Wolfgang L\"uck, and Irakli Patchkoria for stimulating conversations.
We thank Kevin Schreve for explanations about the paper~\cite{Haulmark-Schreve22}.

\section{Preliminaries}
\subsection{Generalised LS-category}

Let~$G$ be a group. A \emph{family~$\F$ of subgroups} of~$G$ is a non-empty set of subgroups of~$G$ that is closed under conjugation and under taking subgroups.
Important examples are the families $\TR$ consisting only of the trivial subgroup, $\FIN$ consisting of all finite subgroups, and $\AME$ consisting of all amenable subgroups. 
For a set~$\calH$ of subgroups of~$G$, the \emph{family~$\F\spann{\calH}$ generated by~$\calH$} is defined as the smallest family containing~$\calH$.
For a family~$\F$ and a subgroup~$H$ of~$G$, we can form the family $\F|_H=\{F\subset H\mid F\in \F\}$ of subgroups of~$H$.

\begin{defn}[{\cite[Definition~2.16]{CLM20}}]\label{defn:generalised LS-cat}
	Let~$X$ be a path-connected space with fundamental group~$G$ and let~$\F$ be a family of subgroups of~$G$.
	A (not necessarily path-connected) open subset~$U$ of~$X$ is an~\emph{$\F$-set} if
	\[
		\im\bigl(\pi_1(U\into X,x)\bigr) \in \F
	\]
	for all~$x\in U$.
	The \emph{generalised LS-category with respect to~$\F$}
	(also \emph{$\F$-category}) $\cat_\F(X)$ is the minimal~$n\in \IN_{\ge 0}$ for which there exists an open cover $X=\bigcup_{i=0}^n U_i$ by $n+1$ many~$\F$-sets.
	If no such finite cover of~$X$ exists, we set~$\cat_\F(X)=\infty$.
	
	The $\F$-category of the group~$G$ is defined as $\cat_\F(G)\coloneqq \cat_\F(BG)$.
\end{defn}

We point out that we use a different normalisation than in~\cite{CLM20}, our value for~$\cat_\F(X)$ is smaller by~$1$.
In the literature similar invariants are sometimes defined in terms of the multiplicity of open covers rather than the cardinality. However, for CW-complexes there is no difference~\cite[Remark~3.13]{CLM20}.  
From here onwards, we will study the generalised LS-category for groups, that is for aspherical spaces (even though some results hold more generally for not necessarily aspherical spaces).

It is a classical result~\cite{Eilenberg-Ganea57,Stallings68,Swan69} that the $\TR$-category $\cat_\TR(G)$ coincides with the cohomological dimension~$\cd(G)$.
The following upper and lower bounds for the~$\F$-category are immediate.

\begin{lem}\label{lem:basic}
	Let~$G$ be a group and let~$\F$ be a family of subgroups of~$G$.
	\begin{enumerate}[label=\enum]
		\item\label{item:anti monotonicity} For a subfamily~$\E\subset \F$, we have $\cat_\F(G)\le \cat_\E(G)$.
		\\
		In particular, $\cat_\F(G)\le \cd(G)$;
		\item\label{item:restriction} For a subgroup~$H\subset G$, we have $\cat_{\F|_H}(H)\le \cat_\F(G)$.
		\\
		In particular, if~$\F|_H=\TR$ then $\cd(H)\le \cat_\F(G)$.
	\end{enumerate}
\end{lem}

Our main object of interest is the $\AME$-category (also \emph{amenable category}) $\cat_\AME(G)$.
A lower bound for the amenable category is given by homology growth.
Recall that a group~$G$ is \emph{of type~$F$} if there exists a finite model for~$BG$.
A group~$G$ is \emph{residually finite} if it admits a \emph{residual chain}~$(\Gamma_i)_{i\in \IN}$, i.e., a nested sequence~$G=\Gamma_0\supset \Gamma_1\supset \Gamma_2\supset\cdots$ such that each~$\Gamma_i$ is a finite index normal subgroup of~$G$ and~$\bigcap_{i\in \IN}\Gamma_i=\{1\}$.
We denote by~$b_k(\Gamma_i;\IF_p)$ the~$k$-th Betti number of~$B\Gamma_i$ with coefficients in~$\IF_p$.

\begin{thm}[{\cite[Theorem~3.2]{Haulmark-Schreve22}\cite[Theorem~1.6]{Sauer16}}]\label{thm:homology growth}
	Let~$G$ be a residually finite group of type~$F$ and let~$(\Gamma_i)_i$ be a residual chain.
	Then we have
	\[
		\lim_{i\to \infty}\frac{b_k(\Gamma_i;\IF_p)}{[G:\Gamma_i]}=0
	\]
	for all~$k>\cat_\AME(G)$ and all primes~$p$.
\end{thm}

\subsection{Classifying spaces for families of subgroups}
Let~$G$ be a group and let~$\F$ be a family of subgroups of~$G$.
A \emph{classifying space}~$\EFG$ for~$G$ with respect to the family~$\F$ is a terminal object in the $G$-homotopy category of $G$-CW-complexes whose isotropy groups lie in~$\F$~\cite{Lueck05_survey}. 
For the trivial family~$\TR$, a model for $E_\TR G$ is given by the universal covering space~$EG$ of~$BG$.
In particular, for every family~$\F$ there is a unique (up to~$G$-homotopy) $G$-map $EG\to \EFG$.

The~$\F$-category of groups can be characterised via classifying spaces for families.
\begin{thm}[{\cite[Proposition~7.5]{CLM20}}]\label{thm:CLM}
	Let~$G$ be a group and let~$\F$ be a family of subgroups of~$G$.
	Then $\cat_\F(G)$ equals the infimum of $n\in \IN_{\ge 0}$ for which the canonical $G$-map $EG\to \EFG$ is $G$-homotopic to a~$G$-map with values in the~$n$-skeleton of~$\EFG$.
\end{thm}

The usual notions of geometric and cohomological dimension of groups admit generalisations to the setting of families.
The \emph{geometric dimension} $\gd_\F(G)$ of~$G$ with respect to~$\F$ is the smallest possible dimension of a model for~$\EFG$. The \emph{cohomological dimension} $\cd_\F(G)$ of~$G$ with respect to~$\F$ is the supremum of degrees in which the~$G$-equivariant Bredon cohomology of~$\EFG$ is non-trivial for some Bredon coefficient module~\cite{Bredon67}.

\begin{cor}\label{cor:upper Bredon}
	Let~$G$ be a group and let~$\F$ be a family of subgroups of~$G$.
	Then we have $\cat_\F(G)\le \gd_\F(G)$.
	
	Moreover, if $\cd_\F(G)=1$ implies $\gd_\F(G)=1$, then $\cat_\F(G)\le \cd_\F(G)$.
	\begin{proof}
		The inequality $\cat_\F(G)\le \gd_\F(G)$ is an immediate consequence of Theorem~\ref{thm:CLM}.
		Since $\gd_\F(G)\le \max\{\cd_\F(G),3\}$~\cite{Lueck-Meintrup00}, it remains to treat the case that $\cd_\F(G)=2$ and $\gd_\F(G)=3$.
		We follow a standard argument using equivariant obstruction theory (see e.g.,~\cite[Theorem~3.6]{GMP19}).
		Let $\EFG$ be a 3-dimensional model and consider the identity map $\id_2\colon (\EFG)_2\to (\EFG)_2$ on its 2-skeleton. The obstruction to extending the restriction $\id_2|_{(\EFG)_1}$ to a $G$-map $\EFG\to (\EFG)_2$ lies in the Bredon cohomology of $\EFG$ in degree~3. This cohomology group is trivial by the assumption that $\cd_\F(G)=2$ and hence there exists a $G$-map $\varphi\colon \EFG\to (\EFG)_2$. 
		By considering the composition
		\[
			EG\to \EFG\xrightarrow{\varphi} (\EFG)_2\into \EFG \,,
		\]
		it follows from Theorem~\ref{thm:CLM} that $\cat_\F(G)\le 2$.
	\end{proof}
\end{cor}
It is conjectured that~$\cd_\F(G)=1$ implies~$\gd_\F(G)=1$ for every family~$\F$, see~\cite{GMP19} for a recent account. While the conjecture is open in general, it is known to hold, e.g., for the family $\FIN$~\cite{Dunwoody79}.

\subsection{Graph products of groups}
Let~$L$ be a flag complex, which shall always mean a finite flag complex, with vertex set~$V$.
Let $G$ be a group and for all~$v\in V$ let $G_v=G$.
The \emph{graph product} $G_L$~\cite{Green90} is the group
\[
	G_L=\ast_{v\in V}G_v/\spann{[G_{v_1},G_{v_2}] \text{ for $v_1,v_2\in V$ spanning an edge in $L$}} \,.
\]
The \emph{right-angled Artin group} (RAAG for short) associated to~$L$ is $A_L=\IZ_L$.
The \emph{right-angled Coxeter group} (RACG for short) associated to~$L$ is $W_L=(\IZ/2\IZ)_L$.

\begin{rem}
	The results of this note hold, when suitably modified, also for graph products with varying vertex groups~$(G_v)_{v\in V}$. However, we restrict ourselves to the case of identical vertex groups for ease of notation.
\end{rem}

For every full subcomplex~$K$ of~$L$, the graph product~$G_L$ retracts onto~$G_K$ by mapping the factors~$(G_v)_v$ corresponding to vertices in~$L\setminus K$ to the trivial element in~$G_K$.

Consider the obvious projection~$q\colon G_L\to \prod_{v\in V}G_v$. 
If~$G$ is abelian, then the kernel of~$q$ is the commutator subgroup~$G_L'$.
Moreover, since the restriction of~$q$ to~$G_\sigma$ is injective for every simplex~$\sigma\subset L$, the intersection of~$G_L'$ with conjugates of~$G_\sigma$ in~$G_L$ is trivial.

By functoriality of the graph product construction $(-)_L$ in the group variable, the projection $\IZ\to \IZ/2\IZ$ induces a map $p\colon A_L\to W_L$ which restricts to the commutator subgroups $p'\colon A_L'\to W_L'$. 
Since $W_L'\subset W_L$ is of finite index and torsion-free, the \emph{virtual cohomological dimension} $\vcd(W_L)$ equals $\cd(W_L')$.

\begin{lem}\label{lem:commutators}
	Let $L$ be a flag complex. The group homomorphism $p'\colon A_L'\to W_L'$ admits a right-inverse.
	In particular, $\vcd(W_L)\le \cd(A_L')$.
	\begin{proof}
		We argue on the level of topological spaces using the polyhedral product construction (see e.g.,~\cite{Panov-Veryovkin16}).
		Models for~$B(A_L')$ and~$B(W_L')$ are given by the polyhedral products~$(\IR,\IZ)^L$ and~$([0,1],\{0,1\})^L$, respectively.
		Consider the map
		\[
			f\colon (\IR,\IZ)\to ([0,1],\{0,1\})\,, \quad x\mapsto\begin{cases}
				x-\lfloor x\rfloor & \text{if $\lfloor x\rfloor$ is even}; \\
				1-(x-\lfloor x\rfloor) & \text{if $\lfloor x\rfloor$ is odd}
			\end{cases}
		\] 
		that ``folds" the real line onto the unit interval.
		The map~$f$ induces a map on polyhedral products and on their fundamental groups the map~$p'\colon A_L'\to W_L'$. 
		A right-inverse to~$f$ is given by the inclusion~$([0,1],\{0,1\})\into (\IR,\IZ)$. It follows that~$W_L'$ is a retract of~$A_L'$ and in particular~$\cd(W_L')\le \cd(A_L')$.
	\end{proof}
\end{lem}

We recall an explicit formula for the virtual cohomological dimension of RACGs \cite[Corollary~8.5.5]{Davis08}.
For the right-angled Coxeter group~$W_L$ associated to a flag complex~$L$, we have
\begin{equation}\label{eqn:vcd}
	\vcd(W_L)=\max\{n\mid \widetilde{H}^{n-1}(L\setminus \sigma;\IZ)\neq 0 \text{ for some simplex }\sigma\subset L \text{ or } \sigma=\emptyset\} \,.
\end{equation}
Here~$\widetilde{H}^*$ denotes reduced cohomology.
The virtual cohomological dimension of RACGs will play a key role due to its interpretation as the cohomological dimension of graph products with respect to the following family.

Let~$G_L$ be a graph product and let~$\FS$ be the family of subgroups of~$G_L$ that is generated by the set of \emph{spherical subgroups} 
\[
	\calS=\{G_\sigma\subset G_L\mid \sigma\subset L \text{ simplex}\} \,.
\]
In the case of RACGs, we have $\FS=\FIN$. In the case of RAAGs, the family~$\FS$ consists of free abelian groups and in particular, $\FS\subset \AME$.

\begin{thm}[{\cite[Corollaries~8.3 and~1.10]{Petrosyan-Prytula20}}]
\label{thm:Petrosyan-Prytula}
	Let~$G_L$ be the graph product associated to a non-trivial group~$G$ and a flag complex~$L$.
	Then~$\cd_\FS(G_L)=\vcd(W_L)$.
	
	Moreover, $\cd_\FS(G_L)=1$ implies~$\gd_\FS(G_L)=1$.
\end{thm}

In view of Theorem~\ref{thm:homology growth}, we recall a computation of homology growth for graph products (using that residually finite amenable groups of type~$F$ are~\emph{$\IF_p$-$\ell^2$-acyclic} by Theorem~\ref{thm:homology growth}). Graph products of residually finite groups are residually finite~\cite{Green90}.

\begin{thm}[{\cite[Theorem~5.1]{Okun-Schreve21}}]
\label{thm:Okun-Schreve}
	Let~$G$ be a non-trivial residually finite amenable group of type~$F$ and 
	let~$G_L$ be the graph product associated to a flag complex~$L$.
	Then, for any residual chain~$(\Gamma_i)_i$ in~$G_L$, we have
	\[
		\lim_{i\to\infty} \frac{b_k(\Gamma_i;\IF_p)}{[G_L:\Gamma_i]} = \widetilde{b}_{k-1}(L;\IF_p) 
	\]
	for all~$k>0$ and all primes~$p$. 
\end{thm}
Here~$\widetilde{b}_{k-1}(L;\IF_p)$ denotes the reduced Betti number of~$L$ with coefficients in~$\IF_p$.
Our formulations of Theorem~\ref{thm:Petrosyan-Prytula} and Theorem~\ref{thm:Okun-Schreve} for graph products are special cases of the results in~\cite{Petrosyan-Prytula20,Okun-Schreve21} which apply to the more general context of group actions with a strict fundamental domain.

\section{Generalised LS-category of right-angled Artin groups}
We investigate the generalised LS-category of RAAGs with respect to several interesting families. Throughout, let $L$ be a finite flag complex.
The $\TR$-category $\cat_\TR(A_L)$ equals the cohomological dimension~$\cd(A_L)$ which is~$\dim(L)+1$. The following lemma provides an upper bound for various families and will be used frequently.

\begin{lem}\label{lem:upper vcd}
	Let~$G_L$ be the graph product associated to a group~$G$ and a flag complex~$L$.
	Let~$\F$ be a family of subgroups of~$G_L$ satisfying~$\FS\subset \F$.
	Then we have $\cat_\F(G_L)\le \vcd(W_L)$.
	\begin{proof}
		Combining Lemma~\ref{lem:basic}~\ref{item:anti monotonicity}, Corollary~\ref{cor:upper Bredon}, and Theorem~\ref{thm:Petrosyan-Prytula} yields the claim.
	\end{proof}
\end{lem}

\subsection{Spherical category}
We compute the generalised LS-category with respect to the family~$\FS$ for RAAGs and RACGs.

\begin{prop}[Spherical category of RAAGs]\label{prop:spherical RAAG}
	Let~$A_L$ be the right-angled Artin group associated to a flag complex~$L$.
	Then we have $\cat_\FS(A_L)=\vcd(W_L)$.
	
	In particular, $\cd(A_L')=\vcd(W_L)$.
	\begin{proof}
		Lemma~\ref{lem:upper vcd} provides the upper bound~$\cat_\FS(A_L)\le \vcd(W_L)$.
		For the lower bound, consider the commutator subgroup~$A_L'$ which satisfies $\FS|_{A_L'}=\TR$. Then Lemma~\ref{lem:basic}~\ref{item:restriction} applied to~$A_L'$ together with Lemma~\ref{lem:commutators} yields
		\[
			\cat_\FS(A_L)\ge \cat_{\FS|_{A_L'}}(A_L')= \cd(A_L') \ge \vcd(W_L) \,,
		\] 
		concluding the proof.
	\end{proof}
\end{prop}
An alternative proof for the lower bound~$\cat_\FS(A_L)\ge \vcd(W_L)$ will be provided by Theorem~\ref{thm:main} below.

A virtually torsion-free group~$G$ satisfies~$\cd_\FIN(G)\ge \vcd(G)$, which follows from the Shapiro lemma for Bredon cohomology. This inequality can be strict, but it is in fact an equality for right-angled Coxeter groups (Theorem~\ref{thm:Petrosyan-Prytula}), as well as for many other examples~\cite{Degrijse-MartinezPerez16}.

\begin{prop}\label{prop:cat fin RACG}
	Let~$G$ be a virtually torsion-free group and suppose~$G$ satisfies~$\cd_\FIN(G)=\vcd(G)$. 
	Then we have $\cat_\FIN(G)=\vcd(G)$.
	\begin{proof}
		By Corollary~\ref{cor:upper Bredon}, we have $\cat_\FIN(G)\le \cd_\FIN(G)=\vcd(G)$.
		The opposite inequality $\vcd(G)\le \cat_\FIN(G)$ follows from Lemma~\ref{lem:basic}~\ref{item:restriction} by restricting to a finite index torsion-free subgroup of~$G$.
	\end{proof}
\end{prop}

\begin{cor}[Finite category of RACGs]
	Let~$W_L$ be the right-angled Coxeter group associated to a flag complex~$L$. Then we have~$\cat_\FIN(W_L)=\vcd(W_L)$.
\end{cor}

\subsection{Amenable category}
We prove the main result of this note. 
\begin{thm}
\label{thm:main}
	Let~$G$ be a non-trivial residually finite amenable group of type~$F$ and let~$G_L$ be the graph product associated to a flag complex~$L$.
	Let~$\F$ be a family of subgroups of~$G_L$ satisfying~$\FS\subset \F\subset \AME$.
	Then we have
	\[
		\cat_\F(G_L)=\vcd(W_L) \,.
	\]
	\begin{proof}
		On the one hand, since~$\FS\subset \F$ we have~$\cat_\F(G_L)\le \vcd(W_L)$ by Lemma~\ref{lem:upper vcd}.
		On the other hand, since~$\F\subset \AME$ we have~$\cat_\F(G_L)\ge \cat_\AME(G_L)$ by Lemma~\ref{lem:basic}~\ref{item:anti monotonicity} and it remains to prove that~$\cat_\AME(G_L)\ge \vcd(W_L)$.
		
		Let~$n=\vcd(W_L)$.
		By~\eqref{eqn:vcd} there exists a (possibly empty) simplex~$\sigma\subset L$ such that~$\widetilde{H}^{n-1}(L\setminus \sigma;\IZ)\neq 0$.
		Let~$K$ be the full subcomplex of~$L$ spanned by the vertices in~$L\setminus \sigma$.
		Then the graph product~$G_K$ associated to~$K$ is a subgroup of~$G_L$.
		Hence~$\cat_\AME(G_L)\ge \cat_\AME(G_K)$ by Lemma~\ref{lem:basic}~\ref{item:restriction} and it suffices to prove $\cat_\AME(G_K)\ge n$.
		
		Indeed, since~$K$ is homotopy equivalent to~$L\setminus \sigma$, we have~$\widetilde{H}^{n-1}(K;\IZ)\neq 0$. 
		The universal coefficient theorem implies that~$\widetilde{H}_{n-1}(K;\IF_p)\neq 0$ for some prime~$p$. 
		By Theorem~\ref{thm:Okun-Schreve} we have
		\[
			\lim_{i\to \infty}\frac{b_n(\Gamma_i;\IF_p)}{[G_K:\Gamma_i]}=\widetilde{b}_{n-1}(K;\IF_p)\neq 0 \,,
		\]
		where~$(\Gamma_i)_i$ is any residual chain in~$G_K$.
		Thus we conclude from Theorem~\ref{thm:homology growth} that $\cat_\AME(G_K)\ge n$. This finishes the proof.
	\end{proof}
\end{thm}

Applying Theorem~\ref{thm:main} to~$G=\IZ$ and~$\F=\AME$ yields the following.

\begin{cor}[Amenable category of RAAGs]
\label{cor:main}
	Let~$A_L$ be the right-angled Artin group associated to a flag complex~$L$. Then we have $\cat_\AME(A_L)=\vcd(W_L)$.
\end{cor}

\subsection{Topological complexity}
Another important generalised LS-category is Farber's topological complexity~\cite{Farber03}.
Let~$G$ be a torsion-free group and consider the product~$G\times G$.
Let~$\F\spann{\Delta}$ be the family of subgroups of~$G\times G$ generated by the diagonal subgroup~$\Delta\subset G\times G$. The \emph{topological complexity} $\TC(G)$ coincides with the $\F\spann{\Delta}$-category~$\cat_{\F\spann{\Delta}}(G\times G)$ of~$G\times G$ by~\cite{FGLO19}, 
which might as well be taken as the definition of~$\TC(G)$.

The topological complexity of RAAGs has been computed~\cite{Cohen-Pruidze08,Gonzalezetal15}.
We recall the precise result for completeness, even though we will not need it in the sequel.
For the right-angled Artin group~$A_L$ associated to a flag complex~$L$ with vertex set~$V$, we have
\[
	\TC(A_L)= \max\{|V_1\cup V_2|\mid V_1,V_2\subset V \text{ each spanning a simplex in }L\} \,.
\]

The topological complexity and amenable category of RAAGs are related by the following inequality, providing positive examples to~\cite[Question~8.1]{CLM20}.

\begin{prop}\label{prop:question}
	Let~$A_L$ be the right-angled Artin group associated to a flag complex~$L$.
	Then we have
	$
		\cat_\AME(A_L\times A_L)\le \TC(A_L)
	$.
	\begin{proof}
		We prove the inequalities $\cat_\AME(A_L\times A_L)\le 2\vcd(W_L)\le \TC(A_L)$.
		Since the product~$A_L\times A_L$ is a right-angled Artin group (associated to the join~$L\ast L$), by Lemma~\ref{lem:upper vcd} we have~$\cat_\AME(A_L\times A_L)\le \vcd(W_L\times W_L)\le 2\vcd(W_L)$.
		
		To prove the remaining inequality, let~$\sigma\subset L$ be a simplex of maximal dimension.
		Consider the subgroup~$A_\sigma\times A_L'$ of~$A_L\times A_L$ which satisfies~$\F\spann{\Delta}|_{A_\sigma\times A_L'}=\TR$. 
		It follows from Lemma~\ref{lem:basic}~\ref{item:restriction} that $\TC(A_L)\ge \cd(A_\sigma\times A_L')$. Since~$A_\sigma$ is free abelian and using Lemma~\ref{lem:commutators}, we obtain
		\[
			\cd(A_\sigma\times A_L')=\cd(A_\sigma)+\cd(A_L') \ge 2\vcd(W_L) \,.
		\]
		This concludes the proof.
	\end{proof}
\end{prop}

\begin{rem}
	The analogous inequality holds for all higher topological complexities~\cite{Rudyak10}, i.e., we have $\cat_\AME((A_L)^r)\le \TC_r(A_L)$ for all~$r\in \IN_{\ge 2}$.
\end{rem}

\subsection{Minimal volume entropy}

The generalised LS-category with respect to families of subgroups with controlled growth is closely related to the (non-)vanishing of minimal volume entropy.
The arguments in this section follow~\cite{Haulmark-Schreve22} to which we refer for the precise definitions.

Let~$G$ be a group of type~$F$ with~$\gd(G)=n$. The \emph{minimal volume entropy}~$\omega(G)$ is defined as the minimal exponential growth rate of balls in cocompact models for~$EG$ of dimension~$n$. There are sufficient conditions for the (non-)vanishing of minimal volume entropy, called the fiber (non-)collapsing assumption (F(N)CA for short)~\cite{Babenko-Sabourau21}. 
More precisely, 
\begin{itemize}
	\item if there exists a finite model for~$BG$ of dimension~$n$ satisfying~FCA, then $\omega(G)=0$;
	\item if every finite model for~$BG$ of dimension~$n$ satisfies~FNCA and~$G$ has uniform uniform exponential growth, then $\omega(G)>0$.
\end{itemize}
The conditions~FCA and~FNCA are not complementary in general.

We will use a reformulation of condition~FCA in the language of generalised LS-category~\cite{Babenko-Sabourau21,Loeh-Moraschini21}.
Let~$G$ be a group as above.
For~$\delta\in \IR_{>0}$, let~$\Subexp_{<\delta}$ denote the family of subgroups~$H$ of~$G$ such that every finitely generated subgroup of~$H$ has subexponential growth with subexponential growth rate~$<\delta$.
It follows from~\cite[Corollary~5.9]{Loeh-Moraschini21} that there exists a finite model for~$BG$ of dimension~$n$ satisfying FCA if and only if
\[
	\cat_{\Subexp_{<(n-k)/n}}(G)< k+1
\]
for some~$k\in \{0,\ldots,n-1\}$. 
(Our values for the generalised LS-category are smaller by~$1$ than in~\cite{Loeh-Moraschini21} because we use a different normalisation.)

The following is a complete characterisation of RAAGs with (non-)vanishing minimal volume entropy. 

\begin{thm}[Minimal volume entropy of RAAGs]
\label{thm:minvolent}
	Let~$A_L$ be the right-angled Artin group associated to a flag complex~$L$ of dimension~$d$. Then~$\widetilde{H}^d(L;\IZ)\neq 0$ if and only if~$\omega(A_L)>0$.
	\begin{proof}
		Using~\eqref{eqn:vcd} we observe that~$\widetilde{H}^d(L;\IZ)\neq 0$ is equivalent to~$\vcd(W_L)=d+1$. 
		For all~$\delta>0$ we have $\FS\subset \Subexp_{<\delta}\subset \AME$ and hence Theorem~\ref{thm:main} implies
		\[
			\cat_{\Subexp_{<\delta}}(A_L)=\vcd(W_L) \,.
		\]
		By the above, we have~$\vcd(W_L)<d+1$ if and only if there exists a finite model for~$B(A_L)$ of dimension~$d+1$ satisfying FCA. In this case~$\omega(A_L)=0$.
		On the other hand, the conditions~FCA and~FNCA are in fact complementary for~$B(A_L)$~\cite{Bregman-Clay21}. Thus, we have~$\vcd(W_L)=d+1$ if and only if every finite model for~$B(A_L)$ of dimension~$d+1$ satisfies~FNCA. In this case~$\omega(A_L)>0$, using that RAAGs have uniform uniform exponential growth. 
	\end{proof}
\end{thm}

Most cases of Theorem~\ref{thm:minvolent} appeared in~\cite[Theorem~1.1]{Haulmark-Schreve22}, which however left open if $\widetilde{H}^d(L;\IZ)=0$ implies $\omega(A_L)=0$ in the case when~$d=2$.
Our resolution of this case goes back to the obstruction theoretical argument in the proof of Corollary~\ref{cor:upper Bredon}.

\section{Graph products of hyperbolic groups}
We provide examples of graph products whose amenable category is maximal, i.e., it equals the cohomological dimension.
For a group~$G$ and a flag complex~$L$ of dimension~$d$,
we have $\cd(G_L)=\cd(G^{d+1})\le (d+1)\cdot \cd(G)$.

\begin{lem}\label{lem:maximality}
	Let~$G$ be a group and let~$L$ be a flag complex of dimension~$d$.
	Let~$\F$ be a family of subgroups of the graph product~$G_L$.
	If there is a simplex~$\sigma\subset L$ of dimension~$d$ such that~$\cat_{\F|_{G_\sigma}}(G_\sigma)=\cd(G_\sigma)$, then $\cat_\F(G_L)=\cd(G^{d+1})$.
	\begin{proof}
		The claim follows at once from Lemma~\ref{lem:basic}~\ref{item:restriction} by restricting to the subgroup~$G_\sigma$ of~$G_L$.
	\end{proof}
\end{lem}

In the following proof we use the notion of simplicial volume of manifolds and some of its standard properties~\cite{Gromov82, Frigerio17}.

\begin{prop}
\label{prop:hyperbolic}
	Let~$G$ be the fundamental group of an oriented closed connected hyperbolic manifold and let~$G_L$ be the graph product associated to a flag complex~$L$.
	Let~$\F$ be a family of subgroups of~$G_L$ satisfying~$\F\subset \AME$.
	Then we have $\cat_\F(G_L)=(\dim(L)+1)\cdot \cd(G)$.
	\begin{proof}
		Let~$d=\dim(L)$ and let~$M$ be an oriented closed connected hyperbolic manifold with~$\pi_1(M)\cong G$.
		Since $\F\subset \AME$, by Lemma~\ref{lem:basic}~\ref{item:anti monotonicity} we have
		\[
			\cat_\AME(G_L)\le \cat_\F(G_L)\le \cd(G_L)\le (d+1)\cdot \cd(G) \,.
		\]
		Since~$M$ is a model for~$BG$, the product~$M^{d+1}$ is a model for~$BG^{d+1}$ and we have~$\cd(G^{d+1})=(d+1)\cdot \cd(G)$.
		Hence by Lemma~\ref{lem:maximality} it suffices to show that $\cat_\AME(M^{d+1})=\dim(M^{d+1})$.
		Indeed, the hyperbolic manifold~$M$ has positive simplicial volume and by the product inequality for simplicial volume so does~$M^{d+1}$. Finally, the amenable category of manifolds with positive simplicial volume is maximal, i.e., it equals the dimension of the manifold.
	\end{proof}
\end{prop}

\begin{rem}
	More generally, Proposition~\ref{prop:hyperbolic} holds by the same proof for every group~$G$ that is the fundamental group of an oriented closed connected aspherical manifold with positive simplicial volume. For a recent list of manifolds known to have positive simplicial volume, see e.g.,~\cite[Example~3.1]{LMR21}.
\end{rem}

\bibliographystyle{alpha}
\bibliography{bib}

\end{document}